\documentstyle[amssymb,righttag,12pt]{amsart}
\newtheorem{theorem}{Theorem}[section]

\newtheorem{corollary} {Corollary}[section]
\newtheorem{lemma} {Lemma}[section]

\newcommand{\re}{{\mathbb R}}

\newcommand{\supp}{\mbox{supp\,}}

\renewcommand{\H}{{\dot{H}}}

\begin{document}
\allowdisplaybreaks

\title[Maximal Operators of Schr\"odinger Type]{Estimates
for multiparameter maximal operators of Schr\"odinger type}

\author{Per Sj\"olin}

 \address{Per Sj\"olin, Department of Mathematics,
Royal Institute of Technology,\hfill\break S-100 44  Stockholm,
Sweden} \email{pers@@math.kth.se}

\author{Fernando Soria}
\address{Fernando Soria,
Departamento de Matem\'aticas, Univ. Aut\'onoma de Madrid,
\hfill\break E-28049 Madrid, Spain} \email{fernando.soria@@uam.es}




\begin{abstract}
Multiparameter maximal estimates are considered for operators of
Schr\"odinger type. Sharp and almost sharp results, that extend  work by Rogers and Villarroya, are obtained. 
 We provide new estimates via the integrability of
the kernel which naturally appears with a $TT^*$ argument and discuss the
behavior at the endpoints. We treat in particular the case of global integrability of the maximal operator on
finite time for solutions to the linear Schr\"odinger equation and make some comments on an open problem 
\end{abstract}

\maketitle

\section{Introduction and main results}

Assuming $a>1$ and letting $f$ belong to the Schwartz class
${\mathcal S}(\re^n)$,
 we set
\begin{equation*}
S_tf(x)=\int_{\re^n} e^{ix\cdot\xi}e^{it|\xi|^a}\widehat
f(\xi)\,d\xi, \qquad x\in \re^n,  \quad t \in \re.
\end{equation*}
Here $\widehat f$ denotes the Fourier transform of the function
$f$, defined by
$$
\widehat f(\xi)=\int_{\re^n} e^{-ix\cdot\xi}f(x)\,dx.
$$

We also set $u(x,t)=(2\pi)^{-n}S_tf(x)$. It then follows that $u(x,0)=f(x)$
and in the case $a=2$, the function $u$ satisfies the Schr\"odinger equation
$i\partial u/\partial t=\Delta u$. 
Also, more generally, if $a=2k$ for some $k=1,2,3,\dots,$ then $u$ satisfies the equation 
$i\partial u/\partial t=\Delta^k u$, if $k$ is odd and $i\partial u/\partial t=-\Delta^k u$, if $k$ is even.

\


We shall study the maximal function $S^*f$
defined by
$$
S^*f(x)=\sup_{0<t<1}\left |S_tf(x)\right |, \qquad x\in \re^n,
$$
and define Sobolev spaces $H_s$ by setting
$$
H_s=\{f\in {\mathcal S}': \, \|f\|_{H_s}<\infty\}, \qquad s\in
\re,
$$
where
$$
\|f\|_{H_s}=\left(\int_{\re^n}\left(1+|\xi|^2\right)^s|\widehat
f(\xi)|^2 d\xi\right)^{1/2}.
$$
The homogeneous Sobolev spaces $\H_s$, for $s\in \re$, are defined
by
$$
\H_s=\{f\in {\mathcal S}': \, \|f\|_{\H_s}<\infty\},
$$
where
$$
\|f\|_{\H_s}=\left(\int_{\re^n} |\xi|^{2s}|\widehat
f(\xi)|^2d\xi\right)^{1/2}.
$$

\

\

The inequality
\begin{equation}
\|S^*f\|_{L^2(B)}\leq C \, \|f\|_{H_s}, \label{local}
\end{equation}
for arbitrary balls $B$ has been studied by several authors. In
the case $n=1$, it is known that (\ref{local}) holds if and only if
$s\geq 1/4$ (see Carleson \cite{C}, Dahlberg and Kenig \cite{DK},
and Sj\"olin \cite{Sj1}). In the case $n=2$ and $a=2$, Lee
\cite{Le}, extending previous results in \cite{TV} and \cite{Ta},
has proved that (\ref{local}) holds for $s>3/8$. In the case
$n\geq 3$, Sj\"olin \cite{Sj1} and Vega \cite{Ve} proved that
(\ref{local}) holds for $s>1/2$. 

\vskip 3mm

As is well known, the inequality (\ref{local}) implies that
$$
\lim_{t\rightarrow 0} \frac 1{(2\pi)^n} \, S_tf(x)=f(x), \quad a.e.,
$$
for every $f\in H_s$. The above estimates therefore give pointwise convergence results. In the case $a=2$, Bourgain \cite{Bou}
has recently improved these results and proved that one has convergence almost everywhere for every $f\in H_s({\mathbb R} ^n)$ if $s>1/2-1/4n$. On the other hand Bourgain has also proved that one does not have convergence almost everywhere for all $f\in H_s({\mathbb R} ^n)$ if $n\ge 5$ and $s<1/2-1/n$. 

\

 For $n=1$
and $a>1$ we set $M^*f=S^*f$  and
$$
M^{**}f(x)=\sup_{t\in \re}\left |S_tf(x)\right |, \qquad x\in
\re.
$$
In harmonic analysis considerable attention has been given to multiparameter singular integrals and related operators. Some examples of this can be seen in the work of E.M. Stein and R. Fefferman \cite{F1}, \cite{F2}, \cite{F3}, \cite{FS}. In this paper we introduce in the same spirit multiparameter operators of Schr\"odinger type. 

\

For $n\geq 2$ and a multiindex $a=(a_1,a_2,\dots,a_n)$, with
$a_j>1$, and $f\in {\mathcal S}(\re^n)$, we now set
$$
S_tf(x)=\int_{\re^n}
e^{ix\cdot\xi}e^{i\left(t_1|\xi_1|^{a_1}+t_2|\xi_2|^{a_2}+\dots+t_n|\xi_n|^{a_n}
\right)}\widehat f(\xi)\,d\xi, \qquad x\in \re^n,
$$
where $t=(t_1,t_2,\dots,t_n) \in \re^n$. In the remaining part of
this paper, $S_t$ will be defined in this way if $n\geq 2$.
Finally, we will define  maximal operators for $n\geq 2$ by
letting
$$
M^{*}f(x)=\sup_{0<t_i<1}\left |S_tf(x)\right |, \qquad x\in \re^n,
$$
and
$$
M^{**}f(x)=\sup_{t_i\in \re}\left |S_tf(x)\right |, \qquad x\in
\re^n.
$$

\

\

In this paper we will study the inequality
\begin{equation}
\|M^{**}f\|_q \leq C \, \|f\|_{\H_s}, \label{two}
\end{equation}
as well as
\begin{equation}
\|M^{*}f\|_q \leq C \, \|f\|_{H_s}, \label{one}
\end{equation}
for different values of $s\in\re, \, 1\leq q\leq \infty $ and the
multiindex $a$. Here we shall use the notation $\|\cdot\|_q=\|\cdot\|_{L^q(\re^n)}$.

\

We can state the following results. The first two theorems are concerned with the case $n=1$. Some parts of them are already known, but we bring them here for the sake of completeness. 

\begin{theorem}
Assume $n=1$. Then the inequality (\ref{two}) holds if and only if
$4\leq q<\infty$ and $s=1/2-1/q$. \label{thm1}
\end{theorem}

Theorem \ref{thm1} is proved in \cite{Sj4}, p.135.

\

\begin{theorem}
Assume $n=1$. Then we have: \\
For $\,1\leq q<2 $ (\ref{one}) holds for no $s$.\\
For $\,q=2$ (\ref{one}) holds for $s>a/4$ and does not hold for $s<a/4$.\\
For $\,2<q<4$ (\ref{one}) holds if and only if $s\geq 1/2-a/4+a/q-1/q$. \\
For $\,4\leq q<\infty$ (\ref{one}) if and only if $s\geq 1/2-1/q$. \\
For $\,q=\infty$ (\ref{one}) holds if and only if $s>1/2$.
\label{thm2}
\end{theorem}

The case $q=2$, $s=a/4$ in the above theorem remains open.
Theorem {\ref{thm2} is well known (see \cite{Sj3}) except for the
case $2<q<4$, $s=1/2-a/4+a/q-1/q$, which has been proved for $a=2$
by Rogers and Villarroya \cite{RV} and will be proved for $a\neq
2$ in this paper.

\

We now considerer the situation of several variables, that is, the multiparameter case. 

\begin{theorem}
Assume $n\geq 2$. Then the inequality (\ref{two}) holds if and
only if $4\leq q<\infty$ and $s=n(1/2-1/q)$. \label{thm3}
\end{theorem}

\

For $n\geq 2$ and $a=(a_1,a_2,\dots,a_n)$ we set
$|a|=a_1+a_2+\dots+a_n$.

\begin{theorem}
Assume $n\geq 2$. Then we have: \\
For $\,1\leq q<2 $ (\ref{one}) holds for no $s$.\\
For $\,q=2$ (\ref{one}) holds for $s>|a|/4$ and does not hold for $s<|a|/4$.\\
For $\,2<q<4$ (\ref{one}) holds if and only if $s\geq n/2-|a|/4+|a|/q-n/q$. \\
For $\,4\leq q<\infty$ (\ref{one}) holds if and only if $s\geq n(1/2-1/q)$. \\
For $\,q=\infty$ (\ref{one}) holds if and only if $s>n/2$.
\label{thm4}
\end{theorem}

In the above theorem, the case $q=2$, $ s=|a|/4$ remains open.

\

Now, set $a=(a_1,a_2,\dots,a_n)$ with $a_j>1$, and set $t=(t_1,t_2,\dots,t_n)$ with $0<t_j<1$. Also, for $\xi_j \in \re$, let $e^{i|\xi_j|^{a_j}}$ have Fourier transform $K^j$. It is known that $K^j\in {\mathcal C}^\infty(\re)$. The function $K^j_{t_j}$ defined as
$$
K^j_{t_j}(x_j)=\frac 1{t_j^{a_j}} K^j\left( \frac {x_j}{t_j^{a_j}}\right ), \quad x_j\in \re,
$$
is then the Fourier transform of $e^{it_j|\xi_j|^{a_j}}$. Hence, $e^{it_1|\xi_1|^{a_1}}e^{it_2|\xi_2|^{a_2}}\dots e^{it_n|\xi_n|^{a_n}}$ has Fourier transform
\begin{equation*}
K_t(x)=K^1_{t_1}(x_1)K^2_{t_2}(x_2)\dots K^n_{t_n}(x_n), \quad x \in \re^n,
\end{equation*}
with $K_t\in {\mathcal C}^\infty \cap {\mathcal S}'$.  Invoking the definition of $ S_t$, we have the identity
\begin{equation*}
S_tf(x)=\int_{\re^n} K_t(y) f(x+y) \, dy=K_t*f(x), \quad {\rm for} \, f \in {\mathcal S}(\re^n). 
\end{equation*}
We then set
\begin{equation*}
T_t f(x)= (2\pi)^{-n} K_t*f(x),  
\end{equation*}
 for $ f \in {L^2 }(\re^n)$ with compact support. Using a standard argument relating maximal funcions and pointwise convergence, one can then prove that Theorem \ref{thm4} has the following consequence.
 
 \
 
 \begin{corollary}
 Assume that $f\in H_{n/4}(\re^n)$ and that $f$ has compact support. Then
 \begin{equation*}
 \lim_{t\to 0} T_t f(x)=f(x), 
 \end{equation*}
 for almost every $x\in \re^n$.
 \end{corollary}

\

We remark that in the above theorems one cannot take $a=1$. In fact, the case $a=1$ is more related to the wave equation than to the Schr\"odinger equation.
We also remark that maximal estimates of the above type have been used to study, among other things, nonlinear equations of Schr\"odinger type.


\

In Section 2 we shall state several lemmas. In Section 3 we will give
proof of these lemmas whereas Section 4 will be devoted to the
proof of the above theorems.
In Section 5 finally we shall make several remarks on the
inequality (\ref{one}) in the open case $n=1,\,a=2,\,q=2,$ and
$s=1/2$.

\

\section{Some lemmas on oscillatory integrals}


In this section we will state several lemmas on oscillatory integrals that may have an interest in their own. They will be used in Section 4 to prove the above theorems.

\begin{lemma}
Assume that $a>1$, $1/2\leq s<1$ and $\mu \in {\mathcal
C}_0^\infty(\re)$. Then,
\begin{equation*}
\left|\int_{\re}
e^{ix\cdot\xi}e^{it|\xi|^a}|\xi|^{-s}\mu(\xi/N)\,d\xi \right|\leq
C\frac1{|x|^{1-s}},
\end{equation*}
for $x\in \re\setminus\{0\}$, $t\in \re$, and
$N=1,2,3,\dots$\label{lemma1}
\end{lemma}

Lemma \ref{lemma1} is contained in \cite{Sj6}.

\

\begin{lemma}
Assume that $a>1$, $1/2\leq \alpha \leq a/2$, $-1<d<1$, and $\mu
\in {\mathcal C}_0^\infty(\re)$. Then,
\begin{equation*}
\left|\int_{\re}
\frac{e^{i(d|\xi|^a-x\xi)}}{(1+\xi^2)^{\alpha/2}}\,\mu(\xi/N)\,d\xi
\right|\leq C\frac1{|x|^{\beta}},
\end{equation*}
for $x\in \re\setminus\{0\}$ and $N=1,2,3,\dots$, where
\begin{equation*}
\beta=\frac{\alpha+a/2-1}{a-1}.
\end{equation*}
\label{lemma2}
\end{lemma}

Observe that the definition of $\beta$ implies that $1/2\leq \beta
\leq 1$. Lemma \ref{lemma2} for $a=2$ is essentially due to Rogers
and Villarroya \cite{RV}. For $a=2$ one has $\beta=\alpha$.

\

\begin{lemma}
Assume that $a>1$, $ \alpha = a/2$, $-1<d<1$,  $\mu \in {\mathcal
C}_0^\infty(\re)$, and $\epsilon>0$. Then,
\begin{equation*}
\left|\int_{\re}
\frac{e^{i(d|\xi|^a-x\xi)}}{(1+\xi^2)^{\alpha/2}[\log(2+\xi^2)]^{1+\epsilon}}\,\mu(\xi/N)\,d\xi
\right|\leq C\, K(x),
\end{equation*}
for $x\in \re\setminus\{0\}$ and $N=1,2,3,\dots$, where $K(x)\in
L^1(\re)$. Moreover, there exists a large constant $C_0$ such that
for $|x|\geq C_0$
\[
K(x)\leq C\frac 1{|x|(\log|x|)^{1+\epsilon}},
\]
whereas for $|x|<C_0$ one has
\begin{eqnarray*}
&(i)&\quad K(x)\leq C, \quad if \quad \alpha\geq 1, \quad and
\\
&(ii)&\quad  K(x)\leq C\frac1{|x|^{1-\alpha}}, \quad if \quad
1/2<\alpha<1.
\end{eqnarray*}
 \label{lemma3}
\end{lemma}

\

Lemma \ref{lemma1} will be used in the case $4\leq q< \infty$,
while Lemma \ref{lemma2} will be used for $2<q<4$. Lemma
\ref{lemma3}, finally, will be used in the case $q=2$.

\

\section{Proofs of the Lemmas}

\noindent{\it Proof of Lemma \ref{lemma2}.} We shall use the
following variants of van der Corput's Lemma (see Stein \cite{St},
p.334):

Assume $a<b$ and set $I=[a,b]$. Let $F\in {\mathcal C}^\infty(I)$
be real valued and let $\psi \in {\mathcal C}^\infty(I)$.

(i) Assume that $|F'(x)|\geq \gamma >0$ for $x\in I$ and that $F'$
is monotonic on I. Then
\begin{equation*}
\left|\int_a^b e^{iF(x)}\psi(x)dx\right|\leq C\frac 1\gamma
\left(|\psi(b)|+\int_a^b|\psi'(x)|dx\right),
\end{equation*}

where C does not depend on $F$, $\psi$ and $I$.

\

(ii) Assume that $|F''(x)|\geq \gamma >0$ for $x\in I$. Then
\begin{equation*}
\left|\int_a^b e^{iF(x)}\psi(x)dx\right|\leq C\frac
1{\gamma^{1/2}} \left(|\psi(b)|+\int_a^b|\psi'(x)|dx\right),
\end{equation*}

where C, again, does not depend on $F$, $\psi$ and $I$.

\

In the proof of the lemma we may assume $d>0$. Clearly, it
suffices to estimate
\begin{equation*}
\left|\int_0^\infty
\frac{e^{i(d|\xi|^a-x\xi)}}{(1+\xi^2)^{\alpha/2}}\,\mu(\xi/N)\,d\xi
\right|.
\end{equation*}

Take first $|x|$ large. Set $F(\xi)=d\xi^a-x\xi$. Then
$F'(\xi)=da\xi^{a-1}-x$ and $F''(\xi)=da(a-1)\xi^{a-2}$. We also
set
\[
\rho=\left(|x|/d\right)^{1/(a-1)}; \quad I_1=[0,\delta\rho];
\quad I_2=[\delta\rho,K\rho], \quad I_3=[K\rho,\infty),
\]
where $\delta$ is to be considered small and $K$ large. On $I_2$
we have for a small positive constant $c$
\[
|F''(x)|\geq c d\left(|x|/d\right)^{(a-2)/(a-1)}.
\]
Setting
\[
\psi(\xi)=(1+\xi^2)^{-\alpha/2}\mu(\xi/N),
\]
we have
\begin{equation*}
\max_{I_2}|\psi|+\int_{I_2}|\psi'|d\xi\leq C
\left(|x|/d\right)^{-\alpha/(a-1)},
\end{equation*}
since $\frac d{d\xi} (\mu(\xi/N)|\leq c\frac 1{1+\xi}$, for
$\xi\geq 0$. van der Corput's Lemma then gives
\begin{equation*}
\left|\int_{I_2} e^{iF(\xi)}\psi(\xi)\,d\xi \right|\leq C
d^{-1/2}\left(\frac{|x|}d\right)^{-\frac{a-2}{2(a-1)}}
\left(\frac{|x|}d\right)^{-\frac{\alpha}{a-1}}
\end{equation*}
$$
=C \frac{d^{\frac{\alpha -1/2}{a-1}}}{|x|^{\frac{\alpha
+a/2-1}{a-1}}} \leq C\frac 1{|x|^\beta},
$$
where we have used in the last inequality that
$\frac{\alpha-1/2}{a-1}\geq 0$ and that $0<d<1$.

\

On $I_1$ we have $\xi \leq \delta \left(|x|/d\right)^{1/(a-1)}$.
Hence, $d \xi^{a-1}\leq \delta^{a-1}|x|$. It follows that
$F'(\xi)|\geq c|x|$ on $I_1$. We also have
\begin{equation*}
\max_{I_1}|\psi|+\int_{I_1}|\psi'|d\xi\leq C ,
\end{equation*}
and van der Corput now gives
\begin{equation*}
\left|\int_{I_1} e^{iF(\xi)}\psi(\xi)\,d\xi \right|\leq C  \frac
1{|x|}\leq C \frac 1{|x|^\beta}.
\end{equation*}

\

On $I_3$ we have $\xi \geq K \left(|x|/d\right)^{1/(a-1)}$. Hence,
$d \xi^{a-1}\geq K^{a-1}|x|$, which implies $F'(\xi)|\geq c|x|$.
Invoking van der Corput again we get
\begin{equation*}
\left|\int_{I_3} e^{iF(\xi)}\psi(\xi)\,d\xi \right|\leq C  \frac
1{|x|}\leq C \frac 1{|x|^\beta}.
\end{equation*}

\

We now consider the case of small values of $x$ ($|x|<C_0$). We
shall consider the cases $\alpha>1$, $1/2\leq \alpha<1$ and
$\alpha=1$ separately. The case $\alpha >1$ is trivial since
\begin{equation*}
\left|\int_0^\infty e^{iF(\xi)}\psi(\xi)\,d\xi \right|\leq
\int_0^\infty |\psi(\xi)| d\xi \leq C \leq C \frac 1{|x|^\beta}.
\end{equation*}
For $1/2\leq \alpha<1$ we use the fact that from the mean value Theorem,
$$0<(1+\xi^2)^{\alpha/2}-\xi^\alpha \leq
(\alpha/2)\xi^{2(\frac\alpha 2-1)}\leq \xi^{\alpha-2}$$ and, therefore,
\[
\frac 1{\xi^\alpha}-\frac 1{(1+\xi^2)^{\alpha/2}}= {\mathcal
O}\left(\frac 1{\xi^{\alpha+2}}\right),
\]
as $\xi\longrightarrow \infty$. It follows that
\[
\int_0^\infty \left|\frac 1{\xi^\alpha}-\frac
1{(1+\xi^2)^{\alpha/2}}\right| d\xi<\infty
\]
and
\begin{equation*}
\left|\int_0^\infty e^{iF(\xi)}\psi(\xi)\,d\xi \right|\leq C+
\left|\int_0^\infty e^{iF(\xi)}\xi^{-\alpha}\mu(\xi/N)\,d\xi
\right|\leq C \frac 1{|x|^{1-\alpha}},
\end{equation*}
where we have used Lemma \ref{lemma1} (replacing the integral
over $\re$ with an integral only on $[0,\infty])$. Observing that
$1-\alpha\leq \beta$, this concludes the case $1/2\leq \alpha<1$.

\

For the case $\alpha=1$ we use the argument in the proof of Lemma
\ref{lemma1} in \cite{Sj6}. Here one obtains
\begin{equation*}
\left|\int_0^\infty e^{iF(\xi)}\psi(\xi)\,d\xi \right|\leq  C
\log\left(\frac 1{|x|}\right), \quad 0<|x|\leq 1/2.
\end{equation*}
and
\begin{equation*}
\left|\int_0^\infty e^{iF(\xi)}\psi(\xi)\,d\xi \right|\leq  C,
\quad 1/2< |x|<C_0.
\end{equation*}
 This finishes the proof of Lemma\ref{lemma2}.

 \hfill{\qed}

\

\

\noindent{\it Proof of Lemma \ref{lemma3}.} We shall first
assume that $|x|$ is large. Choose an even function $\phi_0\in
{\mathcal C}^\infty$ such that $\phi_0(\xi)=1$ for $|\xi|\leq 1/2$
and $\phi(\xi)=0$ for $|\xi|\geq 1$. Set
\[
\psi(\xi)=(1+\xi^2)^{-\alpha/2}[\log(2+\xi^2)]^{-1-\epsilon}\mu(\xi/N),
\]
and $\psi_0=\psi \phi_0$ so that $\supp \psi_0\subset [-1,1]$. We
may assume $d>0$. Let $\rho=(|x|/da)^{1/(a-1)}$. Take $C_0$ large
so that $|x|\geq C_0$ implies $\rho\geq 1000$. Also, take
$C_0>a^2$, $K$ large  and assume $|x|\geq C_0$. Choose $\phi_2\in
{\mathcal C}_0^\infty$ so that $\supp \phi_2\subset
[\rho/4,2K\rho]$ and $\phi_2(\xi)=1$ for $\rho/2\leq \xi\leq
K\rho$. We may also assume that $|\phi_2'(\xi)|\leq C\xi^{-1}$ and
$|\phi_2''(\xi)|\leq C\xi^{-2}$ for $\xi>0$. Set
$\phi_3=(1-\phi_2)\chi_{[K\rho,\infty)}$ and
$\phi_1=(1-\phi_2-\phi_0)\chi_{[0,\rho/2]}$.

\

For $j=1,2,3$, define $\phi_{-j}(\xi)=\phi_j(-\xi)$ and
$F(\xi)=d|\xi|^a-x\xi$. We then have
\begin{equation*}
\int_{-\infty}^\infty e^{iF(\xi)}\psi(\xi)\,d\xi =\sum_{j=-3}^{3}
\int_{-\infty}^\infty e^{iF(\xi)}\psi(\xi)\phi_j(\xi)\,d\xi.
\end{equation*}
The estimates for $j=-1,-2,-3$ can be easily deduced from the
cases $j=1,2,3$, respectively. Setting $\psi_j=\psi \phi_j$,
$j=1,2,3$, we will only consider the integrals
\begin{equation*}
J_j=\int e^{iF(\xi)}\psi_j(\xi)\,d\xi, \quad j=0,1,2,3 .
\end{equation*}
Integrating by parts twice, we get
\begin{equation*}
J_0=\int e^{-ix\xi}e^{id|\xi|^a}\psi_0(\xi)\,d\xi=
\frac{-1}{x^2}\int_{-1}^1 e^{-ix\xi}L(\xi)\,d\xi,
\end{equation*}
where
\begin{equation*}
L(\xi)=\left(\frac{d}{d\xi}\right)^2 (e^{id|\xi|^a}\psi_0(\xi)),
\quad \xi\neq 0.
\end{equation*}
The second integration by parts is justified since
\begin{equation*}
\frac{d}{d\xi} (e^{id|\xi|^a}\psi_0(\xi))=\psi_0'(\xi)e^{id|\xi|^a} +
\psi_0(\xi)ida \,{\rm sign}(\xi) |\xi|^{a-1}e^{id|\xi|^a} \equiv A(\xi)+B(\xi),
\end{equation*}
$A$ and $B$ are both continuous, $A$ is differentiable $\forall \xi$ and $B$ is differentiable  for all $\xi\neq 0$. Moreover, $B(0)=0$ and $B'(\xi)$ is integrable in  $[-1,0)$ and $(0,1]$. We deduce then that
\begin{equation*}
J_0={\mathcal O}\left(\frac 1{x^2}\right),
\end{equation*}
since, for $-1\leq \xi \leq 1$,
\begin{equation*}
L(\xi)= {\mathcal O}\left({|\xi|^{a-2}}\right),
\end{equation*}
and this says, as for $B'$, that $L$  is integrable in $[-1,1]$ when $a>1$.

\

For the remaining estimates we observe that for $j=1,2,3$ and $\xi\geq 1/2$
\[
|\psi_j(\xi)|\leq C (1+\xi^2)^{-\alpha/2}[\log(2+\xi^2)]^{-1-\epsilon}, \quad
\]
\[
|\psi'_j(\xi)|\leq C \xi^{-1}(1+\xi^2)^{-\alpha/2}[\log(2+\xi^2)]^{-1-\epsilon}, \quad
\]
and
\[
|\psi''_j(\xi)|\leq C \xi^{-2}(1+\xi^2)^{-\alpha/2}[\log(2+\xi^2)]^{-1-\epsilon}.
\]
On the interval $[\rho/4, 2K\rho]$ we have
\[
F''(\xi)=da(a-1)\xi^{a-2}\geq
cd\left(\frac{|x|}{d}\right)^{(a-2)/(a-1)},
\]
for a small constant $c>0$. Also,
\[
\max |\psi_2|+\int|\psi'_2|d\xi\leq C \frac 1{\rho^\alpha(\log \rho)^{1+\epsilon}}\leq
C \left(\frac{|x|}{d}\right)^{-\alpha/(a-1)}  \frac 1{(\log |x|)^{1+\epsilon}}.
\]
Using van der Corput's Lemma with the second derivative we obtain
\[
|J_2|\leq C   d^{-1/2}\left(\frac{|x|}{d}\right)^{-(a-2)/2(a-1)}
 \left(\frac{|x|}{d}\right)^{-\alpha/(a-1)}  \frac 1{(\log |x|)^{1+\epsilon}}
 \]
 \[
 =C  d^{1/2} \frac 1{|x|(\log |x|)^{1+\epsilon}} \leq C \frac 1{|x|(\log |x|)^{1+\epsilon}}
\]

\

To estimate $J_1$ observe that $\supp \psi_1\subset [1/2, \rho/2]$. On this interval one has $da\xi^{a-1}\leq da (\rho/2)^{a-1}=2^{1-a}|x|$ and $|F'(\xi)|=|da\xi^{a-1} -x|\geq c|x|\geq cd\xi^{a-1}$. It follows that
\[
\frac {|F''(\xi)|}{|F'(\xi)|}\leq \frac 1\xi, \qquad {\rm and} \qquad
\frac {|F'''(\xi)|}{|F'(\xi)|}\leq \frac 1{\xi^2},
\]
for $1/2 \leq \xi\leq \rho/2$. Integrating by parts twice we obtain
\begin{equation}
J_1=\int e^{iF}\psi_1\,d\xi= \int e^{iF}\frac d{d\xi}\left(\frac 1{iF'} \frac d{d\xi}\left(\frac{\psi_1}{iF'}\right)\right)\,d\xi.
\label{four}
\end{equation}
Now,
\begin{eqnarray*}
\left| \frac d{d\xi}\left(\frac 1{iF'} \frac d{d\xi}\left(\frac{\psi_1}{iF'}\right)\right) \right| &\leq&
\frac {|\psi_1|}{|F'|^2}\left( \frac {|F'''|}{|F'|}+3\frac {|F''|^2}{|F'|^2}\right)+\frac {|\psi_1'|}{|F'|^2}\frac {|F''|}{|F'|}+\frac {|\psi_1''|}{|F'|^2}\\
&=&{\mathcal O}\left(\frac 1{|x|^2\xi^{\alpha+2}}\right).
\end{eqnarray*}
Hence
\[
|J_2|\leq C \int_{1/2}^\infty \frac 1{|x|^2\xi^{\alpha+2}}\, d\xi=
{\mathcal O}\left(\frac 1{|x|^2}\right).
\]

\

It remains to estimate $J_3=\int e^{iF}\psi_3\,d\xi$.  Here $\supp \psi_3\subset [K\rho, \infty]$, and on this interval $da\xi^{a-1}\geq K^{a-1}|x|$ and $|F'(\xi)|\geq c|x|$ and $|F'(\xi)|\geq cda\xi^{a-1}$. Using the same argument (\ref{four}) as for $J_1$ we obtain $|J_3|\leq C/|x|^2$.

\

To finish with the proof of Lemma \ref{lemma3} we must consider
the case $|x|<C_0$. As before, the case $\alpha \geq 1$ is trivial
due to the  integrability of the function $\psi$, and we obtain
$K(x)\leq C$. When $1/2<\alpha<1$, the proof of Lemma
\ref{lemma1} in \cite{Sj6} shows directly that we can take
$K(x)=C/|x|^{1-\alpha}$.

\hfill{\qed}

\

\section{Proofs of the theorems}

\noindent {\it Proof of the case $2<q<4$, $s=1/2-a/4+a/q-1/q$, in
Theorem \ref{thm2}}. Set
\begin{equation*}
Sf(x)=\int_{{\mathbb R}}e^{it(x)|\xi|^a} e^{ix
\xi}\widehat{f}(\xi)d\xi, \quad x\in\re,
\end{equation*}
where $t(x)$ is  measurable and $0<t(x)<1$. We want to prove
\begin{equation*}
\|Sf\|_q\leq C \|f\|_{H_s}= \left(\int_{\mathbb R}|\widehat
f(\xi)|^2(1+|\xi|^2)^{s}d\xi\right)^{1/2}.
\end{equation*}
We set $g(\xi)=\widehat f(\xi)(1+\xi^2)^{s/2}$ and
\begin{equation*}
Tg(x)=\int_{{\mathbb R}}e^{it(x)|\xi|^a} e^{ix
\xi}(1+\xi^2)^{-s/2}g(\xi)d\xi.
\end{equation*}
Then $Sf(x)=Tg(x)$ and it is sufficient to prove that
\begin{equation*}
\|Tg\|_q\leq C \|g\|_{2}.
\end{equation*}
For $N=1,2,3,\dots$ we set
\begin{equation*}
T_Ng(x)=\chi_N(x)\int_{{\mathbb R}}e^{it(x)|\xi|^a} e^{ix
\xi}(1+\xi^2)^{-s/2}\rho_N(\xi)g(\xi)d\xi.
\end{equation*}
Here $\chi_N(x)=\chi(x/N)$ and $\rho_N(x)=\rho(x/N)$, where $\chi$
and $\rho$ are two cut-off functions in ${\mathcal C}_0^\infty$ so
that $\chi(x)=\rho(x)=1$ for $|x|\leq 1$ and $\chi(x)=\rho(x)=0$
for $|x|\geq 2$. They are also assumed to be real-valued. It is
sufficient to prove
\begin{equation*}
\|T_Ng\|_q\leq C \|g\|_{2},
\end{equation*}
with constant $C$ independent of $N$. Its adjoint has the form
\begin{equation*}
T^*_Nh(\xi)=(1+\xi^2)^{-s/2}\rho_N(\xi)\int_{{\mathbb
R}}e^{-it(x)|\xi|^a} e^{-ix \xi}\chi_N(x)h(x)dx.
\end{equation*}
The above is equivalent to prove
\begin{equation}
\|T^*_Nh\|_2\leq C \|h\|_{q'}, \quad N=1,2,3,\dots \label{adj}
\end{equation}
We observe that
\begin{equation*}
\|T^*_Nh\|_2^2=\int\int
I_N(x,y)\chi_N(x)\chi_N(y)h(x)\overline{h(y)}dxdy,
\end{equation*}
where
\[
I_N(x,y)=\int(1+\xi^2)^{-s}e^{i(y-x)\xi}e^{i(t(y)-t(x))|\xi|^a}\mu(\xi/N)d\xi,
\]
and $\mu=\rho^2$.

\

The assumptions $2<q<4$ and $s=1/2-a/4+a/q-1/q$ imply that
$1/4<s<a/4$. Setting $\alpha=2s$ we then have $1/2<\alpha<a/2$.
Lemma \ref{lemma2} then yields
$$
|I_N(x,y)|\leq C|x-y|^{-\beta},
$$
where $\beta=(\alpha+a/2-1)/(a-1)$. It follows that $1/2<\beta<1$.
Set $r=1-\beta$. We define a Riesz potential operator $I_r$ by setting
$$
I_rh(x)=\int_{\re} \frac 1{|x-y|^{1-r}}h(y)dy, \quad x\in \re.
$$
It is not difficult to see that $r=\frac 1{q'}-\frac 1q$, that is
$\frac 1q=\frac 1{q'}-r$. It follows that,
\begin{equation*}
\|I_rh\|_q\leq C \|h\|_{q'}.
\end{equation*}
Hence,
\begin{eqnarray*}
\|T^*_Nh\|_2^2&\leq& C \int\int |x-y|^{-\beta}|h(x)||h(y)|dxdy= C
\int |h(x)|I_r(|h|)(x) dx\\
&\leq& C \|h\|_{q'}\|I_r(|h|)\|_q\leq C \|h\|_{q'}^2.
\end{eqnarray*}
Hence (\ref{adj}) follows and the proof is complete.

\hfill{\qed}

\  

\noindent{\it Proof of Theorem \ref{thm3}}. We first assume
$4\leq q<\infty$ and $s=n(1/2-1/q)$, that is $n/4\leq s<n/2$ and
$q=2n/(n-2s)$. We set
$$
Sf(x)=\int_{\re^n}
e^{ix\cdot\xi}e^{i\left(t_1(x)|\xi_1|^{a_1}+t_2(x)|\xi_2|^{a_2}+\dots+t_n(x)|\xi_n|^{a_n}
\right)}\widehat f(\xi)\,d\xi, \qquad x\in \re^n,
$$
where $t_i(x)$ are measurable and $t_i(x)\in \re$. We want to prove that
\begin{equation*}
\|Sf\|_q \leq C \, \|f\|_{\H_s}=C\left(\int_{\re^n}|\widehat f(\xi)|^2|\xi|^{2s}d\xi\right)^{1/2}.
\end{equation*}
This will follow obviously from the inequality
\begin{equation*}
\|Sf\|_q \leq C \, \left(\int_{\re^n}|\widehat f(\xi)|^2|\xi_1|^{2s/n}\dots
 |\xi_n|^{2s/n}d\xi\right)^{1/2}.
\end{equation*}
Set $g(\xi)=\widehat f(\xi)|\xi_1|^{s/n}\dots  |\xi_n|^{s/n}$ and
\begin{equation*}
Tg(x)= \int_{\re^n}
e^{ix\cdot\xi}e^{i\left(t_1(x)|\xi_1|^{a_1}+t_2(x)|\xi_2|^{a_2}+\dots+t_n(x)|\xi_n|^{a_n}
\right)}g(\xi)|\xi_1|^{-s/n}\dots  |\xi_n|^{-s/n}\,d\xi.
\end{equation*}
Then $Sf(x)=Tg(x)$ and the estimate we want now is
$$
\|Tg\|_q\leq C \|g\|_2.
$$
This can be proved by using Lemma \ref{lemma1} and applying the argument in \cite{Sj6}. We omit the details.

\

We shall now study the necessity of the conditions in Theorem \ref{thm3}. Assume that
$$
\|M^{**}f\|_q\leq C\|f\|_{\H_s}.
$$
Set $f_R(x)=f(Rx)$ for $f\in {\mathcal S}(\re^n)$ and $R>0$. Then $\widehat{f_R}(\xi)=R^{-n}\widehat f(\xi/R)$ and setting $\xi=R\eta$ we obtain for $t=(t_1,t_2,\dots,t_n)$
\begin{eqnarray*}
S_tf_R(x) &=& \int_{\re^n}
e^{ix\cdot\xi}e^{i\left(t_1|\xi_1|^{a_1}+t_2|\xi_2|^{a_2}+\dots+t_n|\xi_n|^{a_n}
\right)}R^{-n}\widehat f(\xi/R)\,d\xi \\
&=& \int_{\re^n}
e^{ix\cdot R\eta}e^{i\left(t_1R^{a_1}|\eta_1|^{a_1}+t_2R^{a_2}|\eta_2|^{a_2}+\dots+t_n
R^{a_n}|\eta_n|^{a_n}
\right)}\widehat f(\eta)\,d\xi \\
&=& S_{\overline t}f(Rx),
\end{eqnarray*}
where ${\overline t}=(t_1R^{a_1},t_2R^{a_2},\dots,t_nR^{a_n})$ It follows that $M^{**}f_R(x)=M^{**}f(Rx)$.
As in \cite{Sj6} one then proves that $q=2n/(n-2s)$ and $s\leq n/2$. A counter-example in \cite{Sj5}, pp. 400-401, shows that the case $s=n/2$ is not possible.

\

It remains to prove that $s\geq n/4$. We shall use a counter-example in \cite{Sj1}, pp. 712-713. Choose $g\in {\mathcal C}_0^\infty(\re)$ with $\int g(\xi)d\xi \neq 0$ and $\supp g \subset [-1,1]$. Define a function $f_v$ for $0<v<1/2$ by the formula
\[
\widehat{f_v}(\xi)=vg(v\xi+1/v), \quad \xi \in \re.
\]
In \cite{Sj1} it is proved that $|S_{t(x)}f_v(x)|\geq c>0$ in a neighbourhood of $x=0$ if $t(x)$ is suitably chosen. Here
\[
S_tf_v(x)=\int_{\re} e^{ix\xi}e^{it|\xi|^a}\widehat{f_v}(\xi)d\xi,
\]
with $a>1$. For $n\geq 2$ we set
\[
f(x)=f_v(x_1)\dots f_v(x_n).
\]
Then
\begin{eqnarray*}
S_tf(x) &=& \int_{\re^n}
e^{ix\cdot\xi}e^{i\left(t_1|\xi_1|^{a_1}+t_2|\xi_2|^{a_2}+\dots+t_n|\xi_n|^{a_n}
\right)}\widehat {f_v}(\xi_1)\dots \widehat {f_v}(\xi_n)\,d\xi \\
&=& S_{t_1}f_v(x_1) \dots S_{t_n}f_v(x_n).
\end{eqnarray*}
Here
\[
S_{t_j}f_v(x_j)=\int_{\re}
e^{ix_j\xi_j}e^{it_j|\xi_j|^{a_j}}\widehat{f_v}(\xi_j)d\xi_j.
\]
It follows that $M^{**}f(x)\geq c>0$ in a neighbourhood of the
origin and hence $\|M^{**}f\|_q\geq c$. On the other hand, it is
easy to see that $\supp \widehat {f_v}$ is included in the
interval $[-1/v^2-1/v,\, -1/v^2+1/v]$. Also,
$|\widehat{f}(\xi)|\leq C v^n$ for all $\xi$. It follows that
\[
\|f\|_{\H_s}^2=\int_{\re^n}|\widehat f(\xi)|^2|\xi|^{2s}d\xi
\leq C v^{-n}v^{2n}v^{-4s}=Cv^{n-4s},
\]
and the right hand side tends to $0$ as $v\rightarrow 0$ if
$n-4s>0$, that is $s<n/4$. Hence, the inequality
$\|M^{**}f\|_q\leq C\|f\|_{\H_s}$ cannot hold for  $s<n/4$. The
proof of Theorem \ref{thm3} is complete.

\hfill{\qed}

\

\

\noindent{\it Proof of Theorem \ref{thm4}}. To treat the case $1\leq q <2$ one can use a counter-example in \cite{Sj3}, pp. 43 and 65. In the case $q=\infty$ we use a counter-example in \cite{Sj5}, pp. 400-401.

The sufficiency in the case $4\leq q<\infty$ follows from Theorem \ref{thm3}. The necessity follows from a counter-example in \cite{Sj3}, pp. 58-59.

We then assume $2<q<4$. We shall prove that inequality (\ref{one}) holds if $s=n/2-|a|/4+|a|/q-n/q$. Set
$$
Sf(x)=\int_{\re^n}
e^{ix\cdot\xi}e^{i\left(t_1(x)|\xi_1|^{a_1}+t_2(x)|\xi_2|^{a_2}+\dots+t_n(x)|\xi_n|^{a_n}
\right)}\widehat f(\xi)\,d\xi, \qquad x\in \re^n,
$$
where each $t_i(x)$ is measurable and $0<t_i(x)<1$. We want to prove that
\begin{equation*}
\|Sf\|_q \leq C \, \|f\|_{H_s}=C\left(\int_{\re^n}|\widehat f(\xi)|^2(1+|\xi|^2)^sd\xi\right)^{1/2}.
\end{equation*}
We have
$$
s=n\frac 12-\frac{a_1+\dots+a_n}4+\frac{a_1+\dots+a_n}q-n\frac 1q=s_1\dots+s_n,
$$
where $s_j=  1/2-a_j/4+a_j/q-1/q$, $j=1,2,\dots,n$. It is sufficient to prove that
\begin{equation*}
\|Sf\|_q \leq C \, \left(\int_{\re^n}|\widehat f(\xi)|^2(1+|\xi_1|^2)^{s_1}\cdots (1+|\xi_n|^2)^{s_n}d\xi\right)^{1/2}.
\end{equation*}
Define $T$ and $T_N$ as before. One obtains
\begin{equation*}
\|T^*_Nh\|_2^2=\int\int
I_N(x,y)\chi_N(x)\chi_N(y)h(x)\overline{h(y)}dxdy,
\end{equation*}
where $I_N(x,y)=I^1_N(x,y)\cdots I^n_N(x,y)$, with
\[
I^j_N(x,y)=\int(1+\xi_j^2)^{-s_j}e^{i(y_j-x_j)\xi_j}e^{i(t_j(y)-t_j(x))|\xi_j|^{a_j}}\mu(\xi_j/N)d\xi_j.
\]
Lemma \ref{lemma2} implies
$$
|I^j_N(x,y)|\leq C|x_j-y_j|^{-\beta_j},
$$
where $\beta_j=(\alpha_j+a_j/2-1)/(a_j-1)$ and $\alpha_j=2s_j$. It follows that
\begin{equation*}
\|T^*_Nh\|_2^2 \leq C\, \int_{\re^n} |h(x)| P_nP_{n-1}\dots P_1(|h|)(x) dx,
\end{equation*}
where
$$
P_jf(x_1, \dots,x_n)=\int_{\re} \frac 1{|x_j-y_j|^{\beta_j}} f(x_1,\dots,x_{j-1}, y_j,x_{j+1},\dots,x_n)dy_j.
$$
We have
$$
\left(\int_{\re}|P_jh(x)|^qdx_j\right)^{1/q} \leq C \left(\int_{\re}|h(x)|^{q'}dx_j\right)^{1/q'},
$$
and the proof can be completed as above (see  also \cite{Sj6}, pp. 407-408).

\

\

Then assume $q=2$. We shall prove that (\ref{one}) holds for $s>|a|/4$. It is sufficient to prove that, for 
 $Sf$ defined as above,
\begin{equation*}
\|Sf\|_2 \leq C_\epsilon \left(\int_{\re^n}|\widehat
f(\xi)|^2(1+\xi^2)^{|a|/4}(\log(2+\xi^2))^{n(1+\epsilon)}d\xi\right)^{1/2}.
\end{equation*}
This, in turn, will follow from the estimate
\begin{equation}
\|Sf\|_2 \leq C_\epsilon \left(\int_{\re^n}|\widehat
f(\xi)|^2 { \Pi}_{j=1}^n\left[ (1+\xi_j^2)^{a_j/4}(\log(2+\xi_j^2))^{1+\epsilon}\right] d\xi\right)^{1/2}. \label{log}
\end{equation}
 We define $T,T_N$ and $I_N$ in the same way as before. Lemma \ref{lemma3} then implies
 $$
 |I_N(x,y)|\leq K_1(x_1-y_1)\cdots K_n(x_n-y_n),
 $$
 with $K_j\in L^1(\re)$ for every $j$. Setting $K(x)=K_1(x_1)\cdots K_n(x_n)$
 we obtain
\begin{eqnarray*}
\|T^*_Nh\|_2^2&\leq& \int |h(x)|K*|h|(x) dx\\
&\leq&  \|h\|_{2}\|K*|h|\|_2\leq C \|h\|_{2}^2.
\end{eqnarray*}
Hence, $T$ is bounded on $L^2$ and (\ref{log}) follows.

\

It remains to prove the necessity for $2\leq q<4$. We shall use a counter-example
 in \cite{Sj2}, pp. 112-113. Let $\phi_j \in {\mathcal
 C}_0^\infty(\re)$ with $\supp \phi_j\subset (-1,\,1)$ and define
 \[
\widehat {f_j}(\xi_j)= \phi_j(N^{a_j/2-1}\xi_j+N^{a_j/2}), \quad
\xi_j\in\re, \quad j=1,2,\dots,n.
 \]
 Here N is large. It is easy to see that $\widehat {f_j}$ vanishes
 outside the interval
 $$
 [-N-N^{1-a_j/2}, -N+N^{1-a_j/2}].
 $$
 We also set $f(x)=f_1(x_1)\cdots f_n(x_n)$. Then
 \[
\|f\|^2_{H_s}\leq C\, \Pi_{j=1}^n (N^{1-a_j/2}) N^{2s}=C
N^{n+2s-|a|/2},
 \]
and hence
\[
\|f\|_{H_s}\leq C\, N^{n/2+s-|a|/4}.
\]
In \cite{Sj2} it is proved that $\phi_j$ can be chosen so that
\[
\sup_{0<t_j<1}\left|\int_{\re}
e^{ix_j\xi_j}e^{it_j|\xi_j|^{a_j}}\widehat{f_j}(\xi_j)d\xi_j\right|\geq
c N^{1-a_j/2},
\]
on a set of measure $\geq c N^{a_j-1}$. It follows that
$$
M^*f(x)\geq c\Pi_{j=1}^n (N^{1-a_j/2}) =c N^{n-|a|/2},
$$
on a set of measure $c \Pi_{j=1}^n (N^{a_j-1})=cN^{|a|-n}$. Hence,
$$
\|M^*f\|_q^q \geq c N^{qn-q|a|/2}N^{|a|-n},
$$
that is
$$
\|M^*f\|_q \geq c N^{n-|a|/2}N^{(|a|-n)/q}.
$$
Inequality (\ref{one}) then implies
$$
 N^{n-|a|/2}N^{(|a|-n)/q} \leq C  N^{n/2+s-|a|/4}.
$$
Letting $N\longrightarrow \infty$ we deduce that
$$
n-\frac{|a|}2+\frac{|a|}q-\frac nq \leq \frac n2+s-\frac{|a|}4,
$$
that is
$$
\frac n2-\frac{|a|}4+\frac{|a|}q-\frac{n}q\leq s.
$$
The proof of Theorem \ref{thm4} is complete.

\hfill{\qed}

\


\section{Some remarks about the case $n=1$, $a=2$, $q=2$}

In this section we assume $n=1$ and $a=2$. It follows from the
method in the proof of Theorem \ref{thm4} that in this case one
has
\begin{equation*}
\|M^*f\|_2 \leq C_\epsilon \left(\int_{\re}|\widehat
f(\xi)|^2(1+\xi^2)^{1/2}(\log(2+\xi^2))^{1+\epsilon}d\xi\right)^{1/2},
\end{equation*}
for $\epsilon>0$

\

As we said above, it remains an open question whether the
logarithmic factor can be removed, that is if
\begin{equation*}
\|M^*f\|_2 \leq C \|f\|_{H_{1/2}}
\end{equation*}
holds. To study this problem we set
\begin{equation*}
Sf(x)=\int_{{\mathbb R}}e^{it(x)\xi^2} e^{ix
\xi}\widehat{f}(\xi)d\xi, \quad x\in\re, \quad f\in {\mathcal S},
\end{equation*}
where $t(x)$ is  measurable and $0<t(x)<1$. We are interested in
the inequality
\begin{equation}
\int_{\re}  |Sf(x)|^2dx\leq C\int_{\mathbb R}|\widehat
f(\xi)|^2(1+|\xi|^2)^{1/2}d\xi, \label{onehalf}
\end{equation}

\

According to Theorem \ref{thm1} one has the estimate
\begin{equation*}
\|Sf\|_4 \leq C \|f\|_{\H_{1/4}},
\end{equation*}
and hence
\begin{equation*}
\int_{|x|\leq 2} \left |Sf(x)\right|^2dx\leq C\int_{\mathbb R}
|\widehat f(\xi) |^2|\xi|^{1/2}d\xi
\end{equation*}
To prove (\ref{onehalf}) it is therefore sufficient to prove that
\begin{equation}
\int_{|x|\geq 2}  |Sf(x)|^2dx\leq C\int_{\mathbb R} |\widehat
f(\xi) |^2(1+|\xi|^2)^{1/2}d\xi, \label{localonehalf}
\end{equation}

\

Consider the phase function $\phi_x(\xi)=t(x)\xi^2+x\xi$. Since
$$\phi'_x(\xi)=2t(x)\xi+x,$$
the zone of \lq\lq non-oscillation" for the kernel of
$S$, that is, when $|\phi_x'(\xi)|\leq 1$,
corresponds to
$$
\{\xi:\,\,|2t(x)\xi+x|\leq 1 \}=\left[\frac{-x-1}{2t(x)},
\frac{-x+1}{2t(x)}\right].
$$
It is therefore natural to look at the operator $L$ defined by
\begin{equation*}
Lf(x)=\int_{\frac{x-1}{2t(x)}}^{\frac{x+1}{2t(x)}}\widehat f(\xi)
d\xi, \quad x\in\re, \quad f\in {\mathcal S}.
\end{equation*}
(For convenience, we have replaced $x$ with $-x$.)

\

We will show that (\ref{localonehalf}) holds with $L$ instead of
$S$. In fact, one has the homogeneous estimate

\begin{theorem}
With the previous notation, we have
\begin{equation}
\int_{|x|\geq 2}\left |Lf(x)\right|^2dx\leq C\int_{\mathbb R}
|\widehat f(\xi) |^2|\xi|d\xi \label{theo3}.
\end{equation}
\label{theo5}
\end{theorem}

\

\noindent{\it Proof}. Setting $g(\xi)=\widehat f(\xi)|\xi|^{1/2}$,
we see that (\ref{theo3}) is equivalent to the inequality
\begin{equation*}
\int_{|x|\geq 2}\left |Tg(x)\right|^2dx\leq C\int_{\mathbb R}\left
|g(\xi)\right |^2d\xi,
\end{equation*}
where
\begin{equation*}
Tg(x)=\int_{\frac{x-1}{2t(x)}}^{\frac{x+1}{2t(x)}} g(\xi)
\frac{d\xi}{|\xi|^{1/2}}, \qquad|x|\geq 2.
\end{equation*}

Observe that the kernel of $T$ is
$$
K_T(x,\xi)=\frac{1}{|\xi|^{1/2}}\chi_{\left[\frac{x-1}{2t(x)},
\frac{x+1}{2t(x)}\right]}(\xi)\chi_{\{|x|\geq 2\}}(x)
$$
Therefore, the kernel of $TT^*$ is
\begin{equation*}
K(x,y)=\int K_T(x,\xi) K_T(y,\xi) d\xi=
\int_{\left[\frac{x-1}{2t(x)}, \frac{x+1}{2t(x)}\right]\cap
\left[\frac{y-1}{2t(y)}, \frac{y+1}{2t(y)}\right]}\frac{d\xi}{|\xi|},
\end{equation*}
for $ |x|, |y|\geq 2$. It follows that
\begin{equation*}
K(x,y) \leq 2\min\left[\log\left(\frac{|x|+1}{|x|-1}\right),
\,\log\left(\frac{|y|+1}{|y|-1}\right)\right], \quad |x|, |y|\geq 2.
\end{equation*}
Using that for $u\geq 2$ one has
$$
\log\frac{u+1}{u-1}=\log\left(1+\frac 2{u-1}\right)\leq \frac
2{u-1}\leq \frac cu,
$$
we get the estimate
$$
|K(x,y)|\leq C\min\left(\frac 1{|x|},\,\frac
1{|y|}\right)\sim\frac 1{|x|+|y|}.
$$
Hence,
$$
|TT^*g(x)|\leq C\left( \frac 1{|x|}\int_{|y|\leq |x|}|g(y)|dy+
\int_{|y|\geq |x|}|g(y)|\frac {dy}{|y|}\right).
$$
The two operators on the right hand side are easily seen to be bounded on $L^p$, for
$1<p<\infty$, and so we obtain the above estimates for $T$ and $L$.

\hfill{\qed}

\

\noindent{\bf Remark}. There is an important feature about the
operator $L$ that we want to point out here and that is that on
scale 1 this is pointwise majorized by the above operator $S$. To
be more precise, we claim that if $f$ is a function so that

\

i) supp $\widehat f \subset [a-1/2,a+1/2]\subset \left[\frac{-x-1}{2t(x)},
\frac{-x+1}{2t(x)}\right]$, for some $a$ and some $x$,  and

ii) $\widehat f$ is positive,

\noindent then
\begin{equation}
|Sf(x)|\geq  c \int_{\frac{-x-1}{2t(x)}}^{\frac{-x+1}{2t(x)}}
\widehat f(\xi) d\xi. \label{claim}
\end{equation}
The reason is simply that
$$
|\phi_x(\xi)-\phi_x(a)|\leq 1/2, \quad \forall \xi \in {\supp}
\widehat f,
$$
(remember that the interval $\left[\frac{-x-1}{2t(x)},
\frac{-x+1}{2t(x)}\right]$ is chosen so that $|\phi_x'(\xi)|\leq 1$
there) and so (\ref{claim}) follows easily.

\

The fact is that many counterexamples in the theory come from the
behavior at the \lq\lq non-oscillation" zone of the kernel
defining $S$. Theorem \ref{theo5} shows that if inequality
(\ref{onehalf}) is not true then we should look for a more
elaborated type of counterexamples.

\

\

In the same spirit of  the above remark, we continue by giving
 a simple proof of the following result
mentioned at the introduction

\begin{theorem} If the inequality
\begin{equation}
\|Sf\|_2 \leq C \|f\|_{H_{s}}\label{nec}
\end{equation}
holds for a constant $C$ independent of $f$ and $t(x)$, then we must have
$s\geq 1/2$.
\end{theorem}

\

\noindent{\it Proof}. We take $M$ large and set
$\widehat f=\chi_{[M,M+1]}$. We let $a=M+1/2$ and choose $x$ so that $-M\leq x\leq -2$. Also set $y=-x$ so that $2\leq y \leq M$. We then take $t(x)=y/2a$. It follows that $0<t(x)<1$, $y/2t(x)=a$, and
$$
\frac{-x-1}{2t(x)}=\frac{y}{2t(x)}-\frac 1{2t(x)}\leq a-1/2=M,
$$
$$
\frac{-x+1}{2t(x)}=\frac{y}{2t(x)}+\frac 1{2t(x)}\geq a+1/2=M+1.
$$
Hence, i) and ii) in the above remark are satisfied and
(\ref{claim}) gives
\[
|Sf(x)|\geq c\int_M^{M+1} d\xi=c.
\]
It follows that
\[
\int |Sf(x)|^2dx\geq c\int_{-M}^{-2}dx\geq cM.
\]
On the other hand
\[
\|f\|^2_{H_s}=\int_M^{M+1} (1+\xi^2)^sd\xi\sim M^{2s}.
\]
If (\ref{nec}) holds, we then get
\[
cM\leq CM^{2s},
\]
and if follows that $1\leq 2s$, that is $s\geq 1/2$.

\hfill{\qed}

\

\

Using the main idea in the previous proof one can also prove the following theorem.
\begin{theorem} Define the maximal operator $U^*$ by setting
\begin{equation}
U^*g(x) \sup_{R>1}\int_{Rx}^{Rx+1} |g(y)|\frac{dy}{|y|^s}, \quad
x\geq 2.
\end{equation}
Then, the inequality
\begin{equation}
\int_2^\infty (U^*g(x))^2dx\leq C \int _{\re} |g(x)|^2dx, \label{theo7}
\end{equation}
holds if and only if $s\geq 1/2$.
\end{theorem}

\

\noindent{\it Proof}. Take $M$ large and let
$g=\chi_{[M,M+1]}$. By taking $R=M/x$ we easily see that
$U^*g(x)\geq 1/(M+1)^s$ on the interval $[2,M]$. Thus, if (\ref{theo7}) holds,
then we should have
$$
(M-2) \frac 1{(M+1)^{2s}}\leq C,\quad {\rm as}\quad  M\longrightarrow \infty.
$$
This implies $2s\geq 1$.

The fact that  (\ref{theo7}) holds if $s\geq 1/2$ follows from Theorem \ref{theo5}.

\hfill{\qed}

\

{\bf Acknowledgement}. The second author was partially supported by DGU Grant MTM2010-16518.

\vskip 2truecm

\end{document}